\documentclass[11pt]{article}
\usepackage{amsmath}
\usepackage{amssymb}
\def\R{\mathbf{R}}

\newtheorem{theo}{\hspace*{\parindent}Theorem}
\newtheorem{lemma}{\hspace*{\parindent}Lemma}
\newtheorem{corol}{\hspace*{\parindent}Corollary}
\newtheorem{prop}{\hspace*{\parindent}Proposition}
\newcounter{theremark}
\newcommand{\rem}{\par\refstepcounter{theremark}\textbf{Remark \arabic{theremark}.} }
\title{Log-convexity and log-concavity of hypergeometric-like functions}
\author{D.\,Karp and S.M.\,Sitnik}
\date{}
\begin{document}
\maketitle

\begin{center}
\parbox{12cm}{
\small\textbf{Abstract.} We find sufficient conditions for
 log-convexity and log-concavity for the functions of the
forms $a\mapsto\sum{f_k}(a)_kx^k$,
$a\mapsto\sum{f_k}\Gamma(a+k)x^k$ and
$a\mapsto\sum{f_k}x^k/(a)_k$. The most useful examples of such
functions are generalized hypergeometric functions.  In
particular, we generalize the Tur\'{a}n inequality for the
confluent hypergeometric function recently proved by Barnard,
Gordy and Richards and log-convexity results for the same function
recently proved by Baricz. Besides, we establish a reverse
inequality which complements naturally the inequality of Barnard,
Gordy and Richards. Similar results are established for the Gauss
and the generalized hypergeometric functions. A conjecture about
monotonicity of a quotient of products of confluent hypergeometric
functions is made.}
\end{center}

\bigskip


\bigskip


\bigskip

\paragraph{1. Introduction.} This paper is motivated by some recent
results dealing with log-convexity and log-concavity of
hypergeometric functions as functions of parameters.  More
specifically, Baricz showed in \cite{Baricz1} that the Kummer
function (or the confluent hypergeometric function)
\begin{equation}\label{eq:F1def}
{_1F_1}(a;c;x)=\sum\limits_{k=0}^{\infty}\frac{(a)_k}{(c)_k}\frac{x^k}{k!},
\end{equation}
where $(a)_k = a(a+1)\cdots(a+k- 1)=\Gamma(a+k)/\Gamma(a)$ is
Pochhammer's symbol, is log-convex in $c$ on $(0,\infty)$ for
$a,x>0$ as well as the function $\mu\mapsto{_1F_1}(a+\mu;c+\mu;x)$
on $[0,\infty)$. This implies, in particular, the reverse
Tur\'{a}n type inequality
\begin{equation}\label{eq:Baricz}
{_1F_1}(a;c+1;x)^2\leq{_1F_1}(a;c;x){_1F_1}(a;c+2;x).
\end{equation}
(This sort of inequalities are called ''Tur\'{a}n type`` after
Paul Tur\'{a}n in a 1946 letter to Szeg\H{o} proved the inequality
$[P_n(x)]^2>P_{n-1}(x)P_{n+1}(x)$, $-1<x<1$, for Legendre
polynomials $P_n$, which has a similar look as (\ref{eq:Baricz})
but different nature (see \cite{Szego}).) Baricz's other results
\cite{Baricz2,Baricz3} deal with log-convexity and some more
general comparisons of means for the Bessel functions (expressible
in terms of ${_0F_1}$) and the Gauss function ${_2F_1}$. Many of
his proofs hinge on the additivity of concavity and logarithmic
convexity. This method does not work, however,  for proving
logarithmic concavity since in general it is not additive.

Closely related results were given a bit earlier by Ismail and
Laforgia in \cite{IsLaf}. In particular, they showed that the
determinant
\[
D_n(x)=\begin{vmatrix}
h(a,c,x) & h(a,c+1,x) & \cdots &  h(a,c+n,x)\\
h(a,c+1,x) & h(a,c+2,x) & \cdots &  h(a,c+n+1,x)\\
\vdots & \vdots & & \vdots\\
 h(a,c+n,x)  &  h(a,c+n+1,x) & \cdots & h(a,c+n,x)
\end{vmatrix},
\]
where $h(a,c,x)=\Gamma(c-a){_1F_1}(a;c;x)/\Gamma(c)$, has positive
power series coefficients.  For $n=1$ this leads to an inequality
which is weaker than (\ref{eq:Baricz}).

 Carey and Gordy conjectured in \cite{CG}
that the Tur\'{a}n type inequality
\[
[{_1F_1}(a;c;x)]^2>{_1F_1}(a+1;c;x){_1F_1}(a-1;c;x),
\]
holds for $a>0$, $c>{a+2}$, $x>{0}$.  Using a clever combination
of contiguous relations and telescoping sums Barnard, Gordy and
Richards have recently shown  in \cite{BGR} that this is indeed
true and even more general inequality
\begin{equation}\label{eq:1F1ineq}
[{_1F_1}(a;c;x)]^2\geq {_1F_1}(a+\nu;c;x){_1F_1}(a-\nu;c;x)
\end{equation}
holds for $a>0$, $c>a\geq\nu-1$ and $x\in{\R}$ or $a\geq\nu-1$,
$c>-1$ ($c\neq{0}$), $x>0$, and positive integer $\nu$.  In fact,
the authors show that the difference of the left-hand and the
right-hand sides of (\ref{eq:1F1ineq}) has positive power series
coefficients for $a>0$, $a\geq\nu-1$, $c>-1$.  They also indicate
that a similar result is true for the generalized hypergeometric
function
\begin{equation}\label{eq:pFqdefined}
{_{p}F_q}\left((a_{p});(b_q);x\right)={_{p}F_q}\left(\left.\!\!\begin{array}{c}(a_{p})\\
(b_q)\end{array}\right|x\!\right)
=\sum\limits_{n=0}^{\infty}\frac{(a_1)_n(a_2)_n\cdots(a_{p})_n}{(b_1)_n\cdots(b_q)_nn!}x^n
\end{equation}
if $p\leq{q+1}$, $a_i>b_i$, $i=2,\ldots,p$ and $b_i>0$,
$i=1,\ldots,q$.

In this paper we shall demonstrate that log-convexity and
log-concavity properties of hypergeometric functions in their
upper and lower parameters serve as an illustration of a more
general phenomenon. Namely, we give sufficient conditions for the
function $x\mapsto{f(a+\delta,x)f(b,x)-f(b+\delta,x)f(a,x)}$ to
have positive power series coefficients if $f(a,x)=\sum
f_k(a)_kx^k$, $f(a,x)=\sum f_k\Gamma(a+k)x^k$ or $f(a,x)=\sum
f_kx^k/(a)_k$, where $f_k>0$ for all $k$.  Log-convexity or
log-concavity then follow immediately.  Section~2 of the paper
contains three general theorems for these three types of functions
and  a corollary which includes direct and reverse Tur\'{a}n type
inequalities. Section~3 collects some applications to
hypergeometric functions. In particular, we extend the result
Barnard, Gordy and Richards (\ref{eq:1F1ineq}) to non-integer
positive $\nu$ and complement it with a reverse inequality giving
asymptotically precise lower bound for the quantity
\[
\frac{{_1F_1}(a+\nu;c;x){_1F_1}(a-\nu;c;x)}{[{_1F_1}(a;c;x)]^2}
\]
bounded by $1$ from above according to (\ref{eq:1F1ineq}). We also
extend some results of Baricz and provide sufficient conditions
for log-concavity and log-convexity of the generalized
hypergeometric functions which are less restrictive then the
conditions $a_i>b_i$ for the (less general) Tur\'{a}n type
inequality given in \cite{BGR}.  We use the generalized Stieltjes
transform representation for ${_{q+1}F_{q}}$ from \cite{KS} to
extend this results to negative $x$. One curious corollary of
these results is positivity of certain finite hypergeometric sums
evaluated at $-1$.

We note in passing the the true Tur\'{a}n type inequalities for
the classical orthogonal polynomials \cite{Gasper,Szego} have been
also shown to exemplify a more general phenomenon. Namely, it has
been demonstrated in  \cite{Szwarc} that they are dependent on
certain monotonicity properties of the coefficients of three-term
recurrence relations. See further development in
\cite{BergSzwarc}.

\paragraph{2. General theorems.}  A function $f:[a,b]\mapsto{\R}$
is said to be Wright-convex (strictly Wright-convex) if
$F_h(x)=f(x+h)-f(x)$ is non-decreasing (increasing) on $[a,b-h]$
for any fixed $h>0$. If $F_h(x)$ is non-increasing (decreasing)
then $f$ is Wright-concave (strictly Wright-concave).  This notion
was introduced by Wright in the 1950s and well studied (see
\cite[p.246]{DragomirPearce} and \cite[p.3]{Mitrinovic}). Clearly,
Wright-convexity implies mid-point convexity and it can be shown
(see \cite[p.2]{Mitrinovic}) that convexity implies
Wright-convexity so that by the celebrated result of Jensen (see,
for instance, \cite[Theorem 1.1.4]{NP}) for continuous functions
all three notions (convexity, Wright-convexity and mid-point
convexity) are equivalent.  In general the inclusions
Convex$\,\subset\,$Wright-convex$\,\subset\,$Midpoint convex are
proper.  We only deal here with log-convexity (log-concavity) of
continuous functions so for our purposes we record

\begin{prop}\label{pr:Wright}
Suppose $f:[a,b]\mapsto\R$ is continuous and has constant sign.
Then $f$ is log-convex  \emph{(}strictly log-convex\emph{)} iff
$x\mapsto{f(x+h)/f(x)}$ is non-decreasing
\emph{(}increasing\emph{)} on $[a,b-h]$ for each fixed $h>0$ and
$f$ is log-concave  \emph{(}strictly log-concave\emph{)} iff
$x\mapsto{f(x+h)/f(x)}$ is non-increasing
\emph{(}decreasing\emph{)} on $[a,b-h]$ for each fixed $h>0$.
\end{prop}
Using Proposition~\ref{pr:Wright} we will formulate our results in
terms of more common log-convexity (log-concavity) while if fact
we prove Wright log-convexity (log-concavity).

In what follows in this section the power series expansions are
understood as formal, so that no questions of convergence are
discussed. It is usually clear in specific applications which
variable range should be considered.  We will use the standard
notation $(a)_n=a(a+1)\cdots(a+n-1)$ for the shifted factorial or
the Pochhammer symbol.  The main idea in the proof of the next
theorem belongs to Fedor Nazarov (University of Wisconsin).

\begin{theo}\label{th:numshift} Let
\begin{equation}\label{eq:f1}
f(a,x)=\sum\limits_{n=0}^{\infty}f_n\frac{(a)_n}{n!}x^n,
\end{equation}
where $f_n>0$ \emph{(}and is independent of $a$\emph{)}. Suppose
$b>a>0$, $\delta>0$.  Then the function
\[
\varphi_{a,b,\delta}(x)=f(a+\delta,x)f(b,x)-f(b+\delta,x)f(a,x)=\sum\limits_{m=2}^{\infty}\varphi_mx^m
\]
has positive power series coefficients $\varphi_m>0$ so that
$a\mapsto{f(a,x)}$ is strictly log-concave for $x>0$ if the
sequence $\{f_{n}/f_{n-1}\}_{n=0}^{\infty}$ is decreasing and
negative power series coefficients $\varphi_m<0$ so that
$a\mapsto{f(a,x)}$ is strictly log-convex for $x>0$ if the
sequence $\{f_{n}/f_{n-1}\}_{n=0}^{\infty}$ is increasing.
\end{theo}
\rem Since, clearly
$\varphi_{a,b,\delta}(x)=-\varphi_{b,a,\delta}(x)$, the sign of
$\varphi_m$ is reversed for $a>b>0$.

\textbf{Proof.} By direct multiplication we have
\[
\varphi_m=\sum\limits_{k=0}^{m}f_{k}f_{m-k}\left\{\frac{(a+\delta)_k(b)_{m-k}}{k!(m-k)!}
-\frac{(b+\delta)_k(a)_{m-k}}{k!(m-k)!}\right\}.
\]
This shows, on inspection, that $\varphi_0=\varphi_1=0$, which
explains why summation starts from $m=2$ in the expansion for
$\varphi_{a,b,\delta}(x)$. Further, we can write $\varphi_m$ in
the form
\begin{equation}\label{eq:GmMk}
\varphi_m=\sum\limits_{k=0}^{[m/2]}f_kf_{m-k}M_k
\end{equation}
with
\[
M_k=\left\{\begin{array}{ll}
[(a+\delta)_k(b)_{m-k}+(a+\delta)_{m-k}(b)_{k}-(a)_k(b+\delta)_{m-k}-(a)_{m-k}(b+\delta)_{k}]/[k!(m-k)!],&
k<m/2

\\[5pt]

[(a+\delta)_k(b)_{m-k}-(a)_k(b+\delta)_{m-k}]/[k!(m-k)!],& k=m/2.
\end{array}\right.
\]
Next we see that
\begin{equation}\label{eq:totalMk}
\sum\limits_{k=0}^{[m/2]}M_k=0
\end{equation}
since for $f_n=1$, $n=0,1,\ldots$, we will have
$f(a,x)=(1-x)^{-a}$ by binomial theorem and hence the left
hand-side of (\ref{eq:totalMk}) is the coefficient at $x^m$ in the
power series expansion of
\[
(1-x)^{-a-\delta}(1-x)^{-b}-(1-x)^{-a}(1-x)^{-b-\delta}=0.
\]
We aim to show that the sequence $\{M_k\}_{k=0}^{[m/2]}$ has
exactly one change of sign, namely some number of initial terms
are negative while all further terms are positive. To establish
the claim note that $(a+\delta)_l(b)_l>(a)_l(b+\delta)_l$ for all
$l$ (since $b>a$ and $x\mapsto(x+\gamma)/x$ is decreasing for
positive $x$ and $\gamma$) and hence $M_0<0$.  Now assume that
$M_k\leq{0}$ for some $k\leq{n/2}$, i.e.
\[
\underbrace{(a+\delta)_k(b)_{n-k}}_{=r}+\underbrace{(a+\delta)_{n-k}(b)_{k}}_{=s}\le\underbrace{(a)_k(b+\delta)_{n-
k}}_{=v}+\underbrace{(a)_{n-k}(b+\delta)_{k}}_{=u}
\]
We want to show that the same inequality is true for $k-1$. We
have by inspection $rs\geq{uv}$ and $v\geq{u}$. A short reflection
shows that together with the above inequality ${r+s}\leq{u+v}$
this yields either $v\geq{r}\geq{s}\geq{u}$ or
$v\geq{s}\geq{r}\geq{u}$ (another apparent possibility
$r>{v}\geq{u}>{s}$ is discarded by noting that it implies that
$u'=u/r$, $v'=v/r$, $s'=s/r$ all belong to $(0,1)$ and satisfy
$u'+v'\geq{1+s'}$ and $s'\geq{u'v'}$ so that $u'+v'\geq{1+u'v'}$
which contradicts the elementary inequality $u'+v'<1+u'v'$,
similarly for $s>{v}\geq{u}>{r}$). We need to prove that
\[
M_{k-1}(\delta)=r\frac{b+m-k}{a+\delta+k-1}+s\frac{a+\delta+m-k}{b+k-1}-v\frac{b+\delta+m-k}{a+k-1}
-u\frac{a+m-k}{b+\delta+k-1}\le{0}
\]
For $\delta=0$, we clearly have
$\frac{b+m-k}{a+k-1}\ge\frac{a+m-k}{b+k-1}$, so the desired
inequality is just a combination of $v\geq{r}$ and  $u+v\geq{r+s}$
with positive coefficients. Treating $u$, $v$, $r$, $s$ as
constants and differentiating with respect to $\delta$, we get
\[
M_{k-1}'(\delta)=u\frac{a+m-k}{(b+\delta+k-1)^2}-v\frac{1}{a+k-1}+s\frac{1}{b+k-1}-r\frac{b+m-k}{(a+\delta+k-1)^2}
\]
which is obviously non-positive since $v\geq{s}$ and $r\geq{u}$,
which proves that $M_{k-1}\le{0}$ and hence that
$\{M_k\}_{k=0}^{[m/2]}$ changes sign exactly once.  Now if
$\{f_{n}/f_{n-1}\}_{n=0}^{\infty}$ is decreasing, then for
$k<m-k+1$
\[
\frac{f_k}{f_{k-1}}>\frac{f_{m-k+1}}{f_{m-k}}~\Leftrightarrow~f_{k}f_{m-k}>f_{k-1}f_{m-k+1}
\]
which combined with (\ref{eq:GmMk}) and (\ref{eq:totalMk}) shows
that $\varphi_m>0$.  Similarly, if
$\{f_{n}/f_{n-1}\}_{n=0}^{\infty}$ is increasing, then for
$k<m-k+1$
\[
\frac{f_k}{f_{k-1}}<\frac{f_{m-k+1}}{f_{m-k}}~\Leftrightarrow~f_{k}f_{m-k}<f_{k-1}f_{m-k+1}
\]
and $\varphi_m<0$. $\square$

\begin{theo}\label{th:numgamma}
Let
\[
g(a,x)=\sum\limits_{n=0}^{\infty}g_n\Gamma(a+n)x^n,
\]
where $g_n>0$ \emph{(}and is independent of $a$\emph{)} and
$\Gamma(\cdot)$ is Euler's gamma function. Suppose $b>a>0$,
$\delta>0$. Then the function
\[
\psi_{a,b,\delta}(x)=g(a+\delta,x)g(b,x)-g(b+\delta,x)g(a,x)=\sum\limits_{m=0}^{\infty}\psi_mx^m
\]
has negative power series coefficients $\psi_m<0$ so that the
function $a\mapsto{g(a,x)}$ is strictly log-convex for $x>0$.
\end{theo}
\noindent\textbf{Proof.} Again, by direct multiplication we have
\[
\psi_m=\sum\limits_{k=0}^{m}g_{k}g_{m-k}\left\{\Gamma(a+\delta+k)\Gamma(b+m-k)
-\Gamma(b+\delta+k)\Gamma(a+m-k)\right\}.
\]
Just like in the proof of Theorem~\ref{th:numshift}, we can write
$\psi_m$ in the form
\begin{equation}\label{eq:HmMk}
\psi_m=\sum\limits_{k=0}^{[m/2]}g_kg_{m-k}M_k
\end{equation}
with
\[
M_k=\left\{\begin{array}{ll}
\underbrace{\Gamma(a+\delta+k)\Gamma(b+k)}_{=r}+\underbrace{\Gamma(a+\delta+m-k)\Gamma(b+k)}_{=s}
\\[5pt]
-\underbrace{\Gamma(a+k)\Gamma(b+\delta+m-k)}_{=v}-\underbrace{\Gamma(a+m-k)\Gamma(b+\delta+k)}_{=u},&
k<m/2

\\[20pt]

\Gamma(a+\delta+k)\Gamma(b+m-k)-\Gamma(a+k)\Gamma(b+\delta+m-k),&
k=m/2.
\end{array}\right.
\]
We aim to show that $M_k<0$ for $k=0,1,\ldots,[m/2]$. The basic
fact that we need is that
$x\mapsto\Gamma(x+\alpha)/\Gamma(x+\beta)$ is strictly increasing
for $x>0$ when $\alpha>\beta\geq{0}$.  This immediately implies
that $M_k<0$ for $k=m/2$ on taking $\alpha=\delta+k=\delta+m-k$,
$\beta=k=m-k$. Further for $k<m-k$ we have
\[
r<v~\Leftrightarrow~\frac{\Gamma(a+\delta+k)}{\Gamma(a+k)}<\frac{\Gamma(b+\delta+m-k)}{\Gamma(b+m-k)},
\]
\[
s<v~\Leftrightarrow~\frac{\Gamma(b+k)}{\Gamma(a+k)}<\frac{\Gamma(b+\delta+m-k)}{\Gamma(a+\delta+m-k)},
\]
\[
u<v~\Leftrightarrow~\frac{\Gamma(b+\delta+k)}{\Gamma(a+k)}<\frac{\Gamma(b+\delta+m-k)}{\Gamma(a+m-k)},
\]
\[
rs<uv~\Leftrightarrow~\frac{\Gamma(b+k)\Gamma(a+\delta+m-k)}{\Gamma(a+k)\Gamma(a+m-k)}
<\frac{\Gamma(b+\delta+k)\Gamma(b+\delta+m-k)}{\Gamma(a+\delta+k)\Gamma(b+m-k)}.
\]
Altogether these inequalities imply $r+s<u+v$ $\Leftrightarrow$
$M_k<0$. Indeed, dividing by $v$ we can rewrite $r+s<u+v$ as
$r'+s'<1+u'$, where $r'=r/v\in(0,1)$, $s'=s/v\in(0,1)$,
$u'=u/v\in(0,1)$.  Since $r's'<u'$ from $rs<uv$, the required
inequality follows from the elementary inequality $r'+s'<1+r's'$.
$\square$

\begin{corol}\label{cr:twosided}
Let $f(a,x)$ be given by \emph{(\ref{eq:f1})} with decreasing
sequence $\{f_{n}/f_{n-1}\}$, then for $b>a>0$ and $x>0$
\[
\frac{\Gamma(a+\delta)\Gamma(b)}{\Gamma(b+\delta)\Gamma(a)}<\frac{f(b+\delta,x)f(a,x)}{f(a+\delta,x)f(b,x)}<1.
\]
\end{corol}
\textbf{Proof.} Indeed, since $(a)_k=\Gamma(a+k)/\Gamma(a)$ we can
choose $g$ in Theorem~\ref{th:numgamma} in the form
$g(a,x)=\Gamma(a)f(a,x)$, where $f$ is given by (\ref{eq:f1}). Now
the estimate from above is just Theorem~\ref{th:numshift} while
the estimate from below is Theorem~\ref{th:numgamma}.  $\square$

\rem Choosing $b=a+\delta$ and $\delta=1$ we get the direct and
reverse Tur\'{a}n type inequalities for $f(a,x)$ given by
(\ref{eq:f1}) with decreasing sequence $\{f_{n}/f_{n-1}\}$:
\[
\frac{a}{a+1}<\frac{f(a+2,x)f(a,x)}{f(a+1,x)^2}<1.
\]

\begin{theo}\label{th:denshift} Let
\[
h(a,x)=\sum\limits_{n=0}^{\infty}\frac{h_n}{(a)_n}x^n,
\]
where $h_n>0$ \emph{(}and is independent of $a$\emph{)}. Suppose
$b>a>0$.  Then the function
\[
\lambda_{a,b,\delta}(x)=h(a+\delta,x)h(b,x)-h(b+\delta,x)h(a,x)=\sum\limits_{m=1}^{\infty}\lambda_mx^m
\]
has negative power series coefficients $\lambda_m<0$ so that the
function $a\mapsto{h(a,x)}$ is strictly log-convex for $x>0$.
\end{theo}
\textbf{Proof.} By direct multiplication we have
\[
\lambda_m=\sum\limits_{k=0}^{m}h_{k}h_{m-k}\left\{\frac{1}{(a+\delta)_k(b)_{m-k}}
-\frac{1}{(b+\delta)_k(a)_{m-k}}\right\}.
\]
This shows, on inspection, that $h_0=0$, which explains why
summation starts from $m=1$ in the expansion for
$\lambda_{a,b,\delta}(x)$. Just like in the proof of
Theorem~\ref{th:numshift}, we can write $\lambda_m$ in the form
\begin{equation}
\lambda_m=\sum\limits_{k=0}^{[m/2]}h_kh_{m-k}M_k
\end{equation}
with
\[
M_k=\left\{\begin{array}{ll}
\underbrace{[(a+\delta)_k(b)_{m-k}]^{-1}}_{=r}+\underbrace{[(a+\delta)_{m-k}(b)_{k}]^{-1}}_{=s}
-\underbrace{[(a)_{m-k}(b+\delta)_{k}]^{-1}}_{=v}-\underbrace{[(a)_k(b+\delta)_{m-k}]^{-1}}_{=u},&
k<m/2

\\[20pt]

[(a+\delta)_k(b)_{m-k}]^{-1}-[(a)_k(b+\delta)_{m-k}]^{-1},& k=m/2.
\end{array}\right.
\]
We aim to show that $M_k<0$ for $k=0,1,\ldots,[m/2]$.  First
$M_{k}<0$ for $k=m-k=m/2$, since
\[
\frac{(b+\delta)_k}{(b)_k}<\frac{(a+\delta)_k}{(a)_k}
\]
because $x\mapsto(x+\alpha)/(x+\beta)$, $\alpha>\beta\geq{0}$, is
decreasing for $x>0$.  Further for $k<m-k$ we have the
inequalities
\[
r<v~\Leftrightarrow~\frac{(b+\delta)_k}{(a+\delta)_k}<\frac{(b)_{m-k}}{(a)_{m-k}},
\]
\[
s<v~\Leftrightarrow~\frac{(b+\delta)_k}{(b)_k}<\frac{(a+\delta)_{m-k}}{(a)_{m-k}},
\]
\[
u<v~\Leftrightarrow~\frac{(b+\delta)_{k}}{(a)_{k}}<\frac{(b+\delta)_{m-k}}{(a)_{m-k}},
\]
\[
rs<uv~\Leftrightarrow~\frac{(b+\delta)_{k}(b+\delta)_{m-k}}{(a+\delta)_{k}(a+\delta)_{m-k}}
<\frac{(b)_{k}(b)_{m-k}}{(a)_{k}(a)_{m-k}}.
\]
Altogether these inequalities imply $r+s<u+v$ $\Leftrightarrow$
$M_k<0$ as shown in the proof of Theorem~\ref{th:numgamma}.
~$\square$

\paragraph{3. Applications to hypergeometric functions}

First consider the Kummer function (\ref{eq:F1def}).  For $a,c>0$
it satisfies the conditions of Theorem~\ref{th:numshift} with
$f_n=1/(c)_n$ and conditions of Theorem~\ref{th:denshift} with
$h_n=(a)_n/n!$. Besides, $\Gamma(a){_1F_1}(a;c;x)$ satisfies the
conditions of Theorem~\ref{th:numgamma} with $g_n=1/[(c)_nn!]$.
Clearly, $f_n/f_{n-1}=1/(c+n-1)$ is decreasing and we are also in
the position to apply Corollary~\ref{cr:twosided}. Using the
Kummer transformation ${_1F_1}(a;c;x)=e^{-x}{_1F_1}(c-a;c;-x)$ we
can extend some of the results to negative $x$. We collect the
consequences of the general theorems for the Kummer function in
the following two statements.

\begin{theo}\label{th:1F1}
Suppose  $\delta>0$.  Then

\emph{a)} for $b>a\geq{0}$, $c>0$ the function
\[
x\mapsto{_1F_1}(a+\delta;c;x){_1F_1}(b;c;x)-{_1F_1}(b+\delta;c;x){_1F_1}(a;c;x)
\]
has positive power series coefficients \emph{(}starting with the
coefficient at $x^2$\emph{)}\emph{;}

\emph{b)} the function
$a\mapsto{_1F_1}(a+\delta;c;x)/{_1F_1}(a;c;x)$ is monotone
decreasing on $[0,\infty)$ for fixed $c,x>0$ and on
$(-\infty,c-\delta]$ for fixed $c>0>x$,  so that
$a\mapsto{_1F_1}(a;c;x)$ is log-concave,
\[
{_1F_1}(a+\delta;c;x)^2\ge{_1F_1}(a;c;x){_1F_1}(a+2\delta;c;x),
\]
on $[0,\infty)$ for fixed $c,x>0$ and on $(-\infty,c]$ for fixed
$c>0>x$\emph{;}

\emph{c)} for $b>a>0$ and $c,x>0$
\begin{equation}\label{eq:1F1twosided1}
\frac{\Gamma(a+\delta)\Gamma(b)}{\Gamma(b+\delta)\Gamma(a)}<\frac{{_1F_1}(b+\delta;c;x){_1F_1}(a;c;x)}{{_1F_1}(a+\delta;c;x){_1F_1}(b;c;x)}<1;
\end{equation}
for $a<b<c-\delta$ and $c>0>x$
\begin{equation}\label{eq:1F1twosided2}
\frac{\Gamma(c-a-\delta)\Gamma(c-b)}{\Gamma(c-b-\delta)\Gamma(c-a)}<\frac{{_1F_1}(b+\delta;c;x){_1F_1}(a;c;x)}{{_1F_1}(a+\delta;c;x){_1F_1}(b;c;x)}<1;
\end{equation}
both sides of both inequalities are sharp in the sense that the
upper bound is attained at $x=0$ and the lower bounds in
\emph{(\ref{eq:1F1twosided1})}, \emph{(\ref{eq:1F1twosided2})} are
attained at $x=+\infty$ and $x=-\infty$, respectively\emph{;}

\emph{d)} for $a>b>0$, $c>0$ and integer $m\geq{2}$
\[
{_4F_3}\left(\left.\!\!\begin{array}{c}-m,a,1-c-m,1-am/(a+b)\\
c,1-b-m,-am/(a+b)\end{array}\right|-1\!\right)>{0}
\]
For $b>a>0$ the sign of inequality is reversed.
\end{theo}

\textbf{Proof.}  Most statements follow immediately from
Theorem~\ref{th:numshift}.  We only need to prove d)  and the
parts of b) and  c) pertaining to negative $x$. Take the
log-concavity inequality from b) and apply the Kummer
transformation ${_1F_1}(a;c;x)=e^{-x}{_1F_1}(c-a;c;-x)$ following
the idea from \cite{BGR}. This yields for $c>a+2\delta$ and $x<0$:
\[
{_1F_1}(c-a-\delta;c;-x)^2-{_1F_1}(c-a;c;-x){_1F_1}(c-a-2\delta;c;-x)
\]\[
=e^{-2x}({_1F_1}(a+\delta;c;x)^2-{_1F_1}(a;c;x){_1F_1}(a+2\delta;c;x))\ge{0}.
\]
In a similar fashion (\ref{eq:1F1twosided2}) follows from
(\ref{eq:1F1twosided1}).  The sharpness of (\ref{eq:1F1twosided1})
and (\ref{eq:1F1twosided2}) at infinity is seen from
\cite[Corollary~4.2.3]{AAR}.

Finally, d)  is a restatement of a) for $\delta=1$, since the
${_4F_3}$ from d) is the coefficient at $x^m$ in the Taylor series
expansion of the function
\[
{_1F_1}(a+1;c;x){_1F_1}(b;c;x)-{_1F_1}(b+1;c;x){_1F_1}(a;c;x).
\]
To see this apply the following easily verifiable identities
\[
(m-k)!=(-1)^k\frac{m!}{(-m)_k},~~(c)_{m-k}=\frac{(-1)^k(c)_m}{(1-c-m)_k},
\]
\[
(a+1)_k(b)_{m-k}-(a)_k(b+1)_{m-k}=
-m(b+1)_{m-1}\frac{(-1)^k(a)_k(1-am/(a+b))_k}{(1-b-m)_k(-am/(a+b))_k}.~~\square
\]

Numerical tests suggests that the following enhancement of
Theorem~\ref{th:1F1} c) is true:

\textbf{Conjecture.}  The ratio in the middle of
(\ref{eq:1F1twosided1}) is monotone decreasing from  $(0,\infty)$
onto $(1,A)$, where $A$ is the left-hand side of
(\ref{eq:1F1twosided1}); the ratio in the middle of
(\ref{eq:1F1twosided2}) is monotone increasing from $(-\infty,0)$
onto $(B,1)$, where $B$ is the left-hand side of
(\ref{eq:1F1twosided2}).

\rem Theorem~\ref{th:1F1} a) and b) for integer $\delta$ and
$b=a+\delta$ recovers \cite[Theorem~1,Corrolary~2]{BGR}.  The
lower bounds in Theorem~\ref{th:1F1} c)  are presumably new. Note
also that although we allow any positive $\delta$ in b), for
integer $\delta$ our parameter ranges are slightly more
restrictive then those from \cite{BGR}.

\rem For the ratio of two Kummer functions with different
denominator parameters Bordelon found in \cite[formula
(5)]{Bordelon} the inequality
\[
1>\frac{{_1F_1}(a;c;x)}{{_1F_1}(a;d;x)}>\frac{\Gamma(c)\Gamma(d-a)}{\Gamma(d)\Gamma(c-a)}
\]
valid for $d>c>a>0$, $x<1$.

\begin{theo}\label{th:1F1-den}
Suppose  $\delta>0$.  Then

\emph{a)} for $d>c>0$ and $a>0$ the function
\[
x\mapsto{_1F_1}(a;c+\delta;x){_1F_1}(a;d;x)-{_1F_1}(a;d+\delta;x){_1F_1}(a;c;x)
\]
has negative power series coefficients \emph{(}starting with the
coefficient at $x$\emph{)}\emph{;}

\emph{b)} the function
$c\mapsto{{_1F_1}(a;c+\delta;x)/{_1F_1}(a;c;x)}$ is monotone
increasing \emph{(}$=$ the function $c\mapsto{{_1F_1}(a;c;x)}$ is
log-convex\emph{)} on $(0,\infty)$ for fixed $a,x>0$ or fixed
$a,x<0$\emph{;}

\emph{c)}  the inequality
${_1F_1}(a+\delta;c+\delta;x)^2\leq{_1F_1}(a+2\delta;c+2\delta;x){_1F_1}(a;c;x)$
holds true for $a\geq{c}>0$, $x>0$ and $a\leq{c}>0$, $x\leq{0}$,
so that $\mu\mapsto{_1F_1}(a+\mu;c+\mu;x)$ is log-convex on
$[0,\infty)$ under these restrictions of parameters.

\emph{d)} for $c>d>0$, $a>0$ and integer $m\geq{1}$
\[
{_4F_3}\left(\left.\!\!\begin{array}{c}-m,-d-m,a,1-cm/(c+d)\\
1-a-m,c+1,-cm/(c+d)\end{array}\right|-1\!\right)>{0}
\]
For $d>c>0$ the sign of inequality is reversed.
\end{theo}

\textbf{Proof.} Statements a) and the part of b) for $a,x>0$
follow from Theorem~\ref{th:denshift}.  The claim c) for positive
$x$ and $a\geq{c}>0$ was proved by Baricz in
\cite[Theorem~2]{Baricz1}. By the Kummer transformation this
yields b) for negative $x$ and $a<0$. Similarly, an application of
the Kummer transformation to the part of b) with $x>0$ gives the
part of c) with $x<0$. Finally, d) is a reformulation of a) for
$\delta=1$.~~$\square$

For the Gauss function
\[
{_2F_1}(a,b;c;x)=\sum\limits_{k=0}^{\infty}\frac{(a)_k(b)_k}{(c)_kk!}x^k
\]
two distinct cases present themselves: if $b>c>0$ then the
sequence
\[
\frac{(b)_{k+1}/(c)_{k+1}}{(b)_{k}/(c)_{k}}=\frac{b+k}{c+k},~k=0,1,2,\ldots,
\]
is decreasing, while for $c>b>0$ it is increasing. We can combine
Theorem~\ref{th:numshift} with Euler's and Pfaff's
transformations,
\[
{_2F_1}(a,b;c;x)=(1-x)^{c-a-b}{_2F_1}(c-a,c-b;c;x)
\]\[
=(1-x)^{-a}{_2F_1}(a,c-b;c;x/(x-1))
=(1-x)^{-b}{_2F_1}(c-a,b;c;x/(x-1)),
\]
to get the following assertions.
\begin{theo}\label{th:2F1decrease}
\emph{a)} Suppose $a'>a\geq{0}$ and $b>c>0$, $\delta>0$. Then the
function
\[
x\mapsto{_2F_1}(a+\delta,b;c;x){_2F_1}(a',b;c;x)-{_2F_1}(a'+\delta,b;c;x){_2F_1}(a,b;c;x)
\]
has positive power series coefficients \emph{(}starting with the
coefficient at $x^2$\emph{)}\emph{;}

\emph{b)} the function $a\mapsto{_2F_1}(a,b;c;x)$ is log-concave,
\[
{_2F_1}(a+\delta,b;c;x)^2\ge{_2F_1}(a,b;c;x){_2F_1}(a+2\delta,b;c;x),
\]
on $[0,\infty)$ for $0<x<1$, $b>c>0$ and $x<0$, $c>0>b$ and on
$(-\infty,c]$ for $0<x<1$,  $c>0>b$ and $x<0$, $b>c>0$\emph{;}

\emph{c)}  For $0<x<1$\emph{:}
\begin{equation}\label{eq:2F1twosided1}
\frac{\Gamma(a+\delta)\Gamma(a')}{\Gamma(a)\Gamma(a'+\delta)}<
\frac{{_2F_1}(a'+\delta,b;c;x){_2F_1}(a,b;c;x)}{{_2F_1}(a+\delta,b;c;x){_2F_1}(a',b;c;x)}<1
\end{equation}
if $b>c>0$, $a'>a>0$ and
\begin{equation}\label{eq:2F1twosided2}
\frac{\Gamma(c-a-\delta)\Gamma(c-a')}{\Gamma(c-a)\Gamma(c-a'-\delta)}<
\frac{{_2F_1}(a'+\delta,b;c;x){_2F_1}(a,b;c;x)}{{_2F_1}(a+\delta,b;c;x){_2F_1}(a',b;c;x)}<1
\end{equation}
if $c>0>b$, $c-\delta>a'>a$.  Further for $x<0$
\emph{(\ref{eq:2F1twosided1})} holds true if $c>0>b$, $a'>a>0$ and
\emph{(\ref{eq:2F1twosided2})} holds true if $b>c>0$,
$c-\delta>a'>a$.
\end{theo}

\begin{theo}\label{th:2F1increase}
\emph{a)} Suppose $a'>a>0$ and $c>b>0$, $\delta>0$. Then the
function
\[
x\mapsto{_2F_1}(a+\delta,b;c;x){_2F_1}(a',b;c;x)-{_2F_1}(a'+\delta,b;c;x){_2F_1}(a,b;c;x)
\]
has negative power series coefficients \emph{(}starting with the
coefficient at $x^2$\emph{)}\emph{;}

\emph{b)} the function $a\mapsto{_2F_1}(a,b;c;x)$ is log-convex,
\[
{_2F_1}(a+\delta,b;c;x)^2\le{_2F_1}(a,b;c;x){_2F_1}(a+2\delta,b;c;x),
\]
on $(-\infty,\infty)$ for $-\infty<x<1$, $c>b>0$.
\end{theo}

\begin{theo}\label{th:2F1den}
\emph{a)} Suppose $d>c>0$ and $a,b,\delta>0$.  Then the function
\[
x\mapsto{_2F_1}(a,b;c+\delta;x){_2F_1}(a,b;d;x)-{_2F_1}(a,b;d+\delta;x){_2F_1}(a,b;c;x)
\]
has negative power series coefficients \emph{(}starting with the
coefficient at $x$\emph{)}\emph{;}

\emph{b)} the function $c\mapsto{{_2F_1}(a,b;c;x)}$ is log-convex
on $(0,\infty)$ if  $a,b>0$, $0<x<1$  or $a,x<0$, $b>0$ or
$b,x<0$, $a>0$\emph{;}

\emph{c)} the function $\mu\mapsto{{_2F_1}(a,b+\mu;c+\mu;x)}$ is
log-convex on $[0,\infty)$ if $a>0$, $b>c>0$, $0<x<1$ or $a<0$,
$b<c>0$, $0<x<1$ or $a>0$, $b<c>0$, $x<0$;

\emph{d)} the function $\mu\mapsto{{_2F_1}(a+\mu,b+\mu;c+\mu;x)}$
is log-convex on $[0,\infty)$ if $a<c, b\geq{c}>0$, $x<0$ or
$a<c$, $b<c$, $c>0$, $0<x<1$.

\end{theo}

\rem Theorems~\ref{th:2F1decrease},\ref{th:2F1increase},
\ref{th:2F1den} do not cover the case when ${_2F_1}$ is expressed
in terms of the Jacobi polynomials so that the original Tur\'{a}n
type inequalities due to  Szeg\H{o} \cite{Szego} and Gasper
\cite{Gasper} cannot be derived from it.

Applications to generalized hypergeometric function
${_{q+1}F_{q}}$ hinge on the following observation which might be
of independent interest:

\begin{lemma}\label{lm:general}
Let $A(x)=a_0+a_1x+\cdots+a_nx^n$ and
$B(x)=b_0+b_1x+\cdots+b_nx^n$ have positive coefficients. Then
$A'(x)B(x)-B'(x)A(x)$ has non-negative coefficients if
\begin{equation}\label{eq:chain}
\frac{a_n}{b_n}\geq\frac{a_{n-1}}{b_{n-1}}\geq\cdots\geq\frac{a_1}{b_1}\geq\frac{a_0}{b_0}
\end{equation}
and non-positive coefficients if
\begin{equation}\label{eq:chain1}
\frac{a_n}{b_n}\leq\frac{a_{n-1}}{b_{n-1}}\leq\cdots\leq\frac{a_1}{b_1}\leq\frac{a_0}{b_0}.
\end{equation}
If $A$ and $B$ are not identical then some of the coefficients are
positive under \emph{(\ref{eq:chain})} \emph{(}so that
$A'(x)B(x)-B'(x)A(x)>0$ for $x>0$\emph{)} and negative under
\emph{(\ref{eq:chain1})} \emph{(}so that $A'(x)B(x)-B'(x)A(x)<0$
for $x>0$\emph{)}.
\end{lemma}

\textbf{Proof.}  We have
\[
A'(x)B(x)-B'(x)A(x)=\sum\limits_{k=1}^{n}\sum\limits_{i=0}^{n}ka_{k}b_{i}x^{i+k-1}
-\sum\limits_{k=1}^{n}\sum\limits_{i=0}^{n}kb_{k}a_{i}x^{i+k-1}
\]\[
=\frac{1}{x}\sum\limits_{k=1}^{n}\sum\limits_{i=0}^{n}x^{i+k}k(a_{k}b_{i}-a_{i}b_{k})
=\sum\limits_{m=1}^{2n}x^{m-1}\sum\limits_{\substack{i+k=m\\i\geq{0},k\geq{1}}}^{k,i\leq{n}}k(a_{k}b_{i}-a_{i}b_{k})
\]
Since each term in the inner sum with $i=k$ is clearly zero, we
may write:
\[
\sum\limits_{\substack{i+k=m\\i\geq{0},k\geq{1}}}^{k,i\leq{n}}k(a_{k}b_{i}-a_{i}b_{k})
=\sum\limits_{\substack{i+k=m,
k<i\\i\geq{0},k\geq{1}}}^{k,i\leq{n}}k(a_{k}b_{i}-a_{i}b_{k})+\sum\limits_{\substack{i+k=m,k>i\\i\geq{0},k\geq{1}}}^{k,i\leq{n}}k(a_{k}b_{i}-a_{i}b_{k})
\]
Due to condition (\ref{eq:chain}) every term in the second sum is
non-negative.  For each term in the first sum (say indexed
$k=k^*$, $i=i^*$, $k^*<i^*$), there is a term in the second sum
with $k=i^*$, $i=k^*$ and
\[
k^*(a_{k^*}b_{i^*}-a_{i^*}b_{k^*})+i^*(a_{i^*}b_{k^*}-a_{k^*}b_{i^*})=
(i^*-k^*)(a_{i^*}b_{k^*}-a_{k^*}b_{i^*})\geq{0}
\]
Since $A$ and $B$ are not identical at least one of the
inequalities (\ref{eq:chain}) is strict which implies that some of
the terms above is strictly positive.  The second statement
follows from $A'(x)B(x)-B'(x)A(x)=-(B'(x)A(x)-A'(x)B(x))$ by
exchanging the roles of $A$ and $B$.~~$\square$

\rem\label{rm:n12} It is easy to verify directly that conditions
(\ref{eq:chain}) are also necessary for $A'(x)B(x)-B'(x)A(x)>0$ if
$x>0$ when $n=1$ and $n=2$ (similarly for conditions
(\ref{eq:chain1})).

\rem After this paper was already written \'{A}rp\'{a}d Baricz
brought the recently published article \cite{HVV} to our
attention. Lemma~\ref{lm:general} repeats almost precisely the
contents of Theorem~4.4 from \cite{HVV}.  Interestingly, the proof
in \cite{HVV} is by induction and hence differs substantially from
ours.

Now consider
\[
R(x)=\frac{\prod\limits_{k=1}^{q}(a_k+x)}{\prod\limits_{k=1}^{q}(b_k+x)}
\]
with positive $a_k$, $b_k$. Let $e_m(c_1,\ldots,c_q)$ denote
$m$-th elementary symmetric polynomial,
\[
e_1(c_1,\ldots,c_q)=c_1+c_2+\cdots+c_q,
\]\[
e_2(c_1,\ldots,c_q)=c_1c_2+c_1c_3+\cdots+c_1c_q+c_2c_3+c_2c_4+\cdots+c_2c_q+\cdots+c_{q-1}c_q,
\]\[
\ldots,\]
\[
e_q(c_1,\ldots,c_q)=c_1c_2\cdots{c_q}.
\]

\begin{lemma}\label{lm:ratio}
The function $R(x)$ is monotone increasing on $(0,\infty)$ if
\begin{equation}\label{eq:symmetric-chain}
\frac{e_q(b_1,\ldots,b_q)}{e_q(a_1,\ldots,a_q)}\geq
\frac{e_{q-1}(b_1,\ldots,b_q)}{e_{q-1}(a_1,\ldots,a_q)}\geq\cdots\geq
\frac{e_1(b_1,\ldots,b_q)}{e_1(a_1,\ldots,a_q)}\geq{1}
\end{equation}
and monotone decreasing if
\begin{equation}\label{eq:symmetric-chain1}
\frac{e_q(b_1,\ldots,b_q)}{e_q(a_1,\ldots,a_q)}\leq
\frac{e_{q-1}(b_1,\ldots,b_q)}{e_{q-1}(a_1,\ldots,a_q)}\leq\cdots\leq
\frac{e_1(b_1,\ldots,b_q)}{e_1(a_1,\ldots,a_q)}\leq{1}.
\end{equation}
\end{lemma}
\textbf{Proof.} We have
\[
A(x)\equiv\prod\limits_{k=1}^{q}(a_k+x)=\sum\limits_{k=0}^{q}e_k(a_1,\ldots,a_q)x^{q-k}=\sum\limits_{k=0}^{q}e_{q-k}(a_1,\ldots,a_q)x^{k}
\]
\[
B(x)\equiv\prod\limits_{k=1}^{q}(b_k+x)=\sum\limits_{k=0}^{q}e_k(b_1,\ldots,b_q)x^{q-k}=\sum\limits_{k=0}^{q}e_{q-k}(b_1,\ldots,b_q)x^{k}
\]
\[
R'(x)=\frac{A'(x)B(x)-B'(x)A(x)}{[B(x)]^2}
\]
Hence by Lemma~\ref{lm:general} we can assert that
\[
1=\frac{e_0(a_1,\ldots,a_q)}{e_0(b_1,\ldots,b_q)}\geq\frac{e_1(a_1,\ldots,a_q)}{e_1(b_1,\ldots,b_q)}\geq\frac{e_2(a_1,\ldots,a_q)}{e_2(b_1,\ldots,b_q)}\geq\cdots\geq\frac{e_q(a_1,\ldots,a_q)}{e_q(b_1,\ldots,b_q)}
\]
which is the same as (\ref{eq:symmetric-chain}), is sufficient for
$R'(x)>0$.  Similarly for $R'(x)<0$.~$\square$

\rem Conditions $b_i>a_i$, $i=1,\ldots,q$ are clearly sufficient
but not necessary for (\ref{eq:symmetric-chain}).

\begin{theo}\label{th:pFq}
Put
\[
f(\alpha,x)={_{q+1}F_q}\left(\alpha,(a_{q});(b_q);x\right)
\]
and suppose $a_i,b_i>0$, $i=1,\ldots,q$, $\beta>\alpha>0$. Then

\emph{a)} for any fixed $\delta>0$ the function
\[
x\mapsto f(\alpha+\delta,x)f(\beta,x)-f(\beta+\delta,x)f(\alpha,x)
\]
has negative power series coefficients if
\emph{(\ref{eq:symmetric-chain})} holds and positive power series
coefficients if \emph{(\ref{eq:symmetric-chain1})} holds\emph{;}

\emph{b)} under condition \emph{(\ref{eq:symmetric-chain})} and
$x\in(0,1)$ the function
$\alpha\mapsto{f(\alpha+\delta,x)/f(\alpha,x)}$ is monotone
increasing on $[0,\infty)$ for any fixed $\delta>0$ so that the
function $\alpha\mapsto{f(\alpha,x)}$ is log-convex\emph{;}

\emph{c)} under condition \emph{(\ref{eq:symmetric-chain1})} and
$x\in(0,1)$ the function
$\alpha\mapsto{f(\alpha+\delta,x)/f(\alpha,x)}$ is monotone
decreasing on $[0,\infty)$ for any fixed $\delta>0$ so that the
function $\alpha\mapsto{f(\alpha,x)}$ is log-concave\emph{;}

\emph{d)} under condition \emph{(\ref{eq:symmetric-chain1})} and
$x\in(0,1)$
\[
\frac{\Gamma(\alpha+\delta)\Gamma(\beta)}{\Gamma(\beta+\delta)\Gamma(\alpha)}
<\frac{f(\beta+\delta,x)f(\alpha,x)}{f(\alpha+\delta,x)f(\beta,x)}<1.
\]
\end{theo}

\rem For ${_3F_2}\left(\alpha,a_1,a_2;b_1,b_2;x\right)$ conditions
(\ref{eq:symmetric-chain}) read:
\[
\frac{b_1b_2}{a_1a_2}\geq\frac{b_1+b_2}{a_1+a_2}\geq{1}
\]
and for condition (\ref{eq:symmetric-chain1}) both inequalities
are reversed.  According to remark~\ref{rm:n12} after
Lemma~\ref{lm:general} in this case these conditions are both
necessary and sufficient for the increase or decrease of the
function $x\mapsto(a_1+x)(a_2+x)/[(b_1+x)(b_2+x)]$.

We can extend Theorem~\ref{th:pFq}(b) to negative $x$ using the
generalized Stieltjes transform representation
\begin{equation}\label{eq:Frepr}
{_{q+1}F_q}\left(\left.\!\!\begin{array}{c}\alpha,(a_{q})\\
(b_q)\end{array}\right|x\!\right)
=\!\!\int\limits_{0}^{1}\frac{\rho((a_q);(b_q);s)ds}{(1-sx)^{\alpha}},
\end{equation}
valid for $b_i>a_i>0$, $i=1,\ldots,q$, and $x<1$, recently
obtained by the authors in \cite{KS}. An explicit expression for
the positive function $\rho(s)$ is given in \cite{KS}. Indeed, the
inequality
\[
f(\alpha+\delta,x)f(\beta,x)<f(\beta+\delta,x)f(\alpha,x)
\]
is exactly the Chebyshev inequality \cite[Chapter IX, formula
(1.1)]{Mitrinovic}.  Hence, Theorem~\ref{th:pFq}(b) is true for
all $x<1$ if $b_i>a_i>0$, $i=1,\ldots,q$.

The claim of Lemma~\ref{lm:general} made under condition
(\ref{eq:chain1}) clearly remains true when
$a_n=a_{n-1}=a_{n-r}=0$. Condition (\ref{eq:symmetric-chain1})
then reads
\[
\frac{e_q(b_1,\ldots,b_q)}{e_{q-r-1}(a_1,\ldots,a_{q-r-1})}\leq
\frac{e_{q-1}(b_1,\ldots,b_q)}{e_{q-r-2}(a_1,\ldots,a_{q-r-1})}\leq\cdots
\leq\frac{e_{r+2}(b_1,\ldots,b_q)}{e_{1}(a_1,\ldots,a_{q-r-1})}\leq
e_{r+1}(b_1,\ldots,b_q)
\]
This permits the application of the Theorem~\ref{th:numshift} to
${_pF_q}$ with $p\leq{q}$.  Of course in this case we will only
have the decreasing sequence $\{F_n/F_{n-1}\}$ and log-concavity
of ${_pF_q}$ in the upper parameters. For instance, an analogue of
Theorem~\ref{th:1F1}(d) for ${_qF_q}$ is the inequality
\[
{_{2q+2}F_{2q+1}}\left(\left.\!\!\begin{array}{c}-m,\alpha,a_1,\ldots,a_{q-1},1-b_1-m,\ldots,1-b_q-m,1-\alpha{m}/(\alpha+\beta)\\
b_1,\ldots,b_q,1-a_1-m,\ldots,1-a_{q-1}-m,1-\beta-m,-\alpha{m}/(\alpha+\beta)\end{array}\right|-1\!\right)>{0}
\]
valid for $\alpha>\beta>0$, integer $q\geq{1}$, integer $m\geq{2}$
and
\[
\frac{e_q(b_1,\ldots,b_q)}{e_{q-1}(a_1,\ldots,a_{q-1})}\leq
\frac{e_{q-1}(b_1,\ldots,b_q)}{e_{q-2}(a_1,\ldots,a_{q-1})}\leq\cdots
\leq\frac{e_{2}(b_1,\ldots,b_q)}{e_{1}(a_1,\ldots,a_{q-1})}\leq
e_{1}(b_1,\ldots,b_q),
\]
where all $a_i,b_i>0$.  For ${_2F_2}(\alpha,a_1;b_1,b_2;x)$ these
inequalities simplify to the single condition
$a_1\geq{b_1b_2/(b_1+b_2)}$ which ensures log-concavity in
$\alpha$.
 Theorems~\ref{th:numgamma} and \ref{th:denshift} can be applied as
well.

\rem Note that some results related to log-convexity of the
modified Struve functions expressible in terms of ${_1F_2}$ are
given in \cite{JoshiNalwaya}.

\paragraph{5. Acknowledgements.} We thank Fedor Nazarov whose idea
led to the proof of Theorem~\ref{th:numshift}, \'{A}rp\'{a}d
Baricz for providing us with the copies of
\cite{Baricz1,Baricz2,Baricz3,HVV,Bordelon} and many valuable
comments, Matti Vuorinen for sharing the problem of log-concavity
of ${_1F_1}$ with us and Michael Gordy for sending us the copies
of \cite{CG,BGR} and other useful communication. The work of the
first author has been supported by Far Eastern Branch of the
Russian Academy of Sciences (grants 09-III-A-01-008 and
09-II-CO-01-003), Russian Basic Research Fund (grant
08-01-00028-a) and the Presidential Grant for Leading Scientific
Schools (grant 2810.2008.1).

\end{document}